\documentclass{amsart}
\usepackage{amsmath, amscd, amssymb, latexsym, hhline, euscript, mathrsfs, bbm}
\usepackage[all]{xy}
\usepackage{enumitem}
\usepackage[usenames]{color} 

\def\1{{\mathbbm{1}}}

\def\w{{\omega}}

\def\CC{\mathcal{C}}

\def\Hom{{\mbox{\rm Hom}}}
\def\Aut{{\mbox{\rm Aut}}}

\def\Vec{{\mbox{\rm Vec}}}

\theoremstyle{definition}

\theoremstyle{remark}

\makeatletter
\def\namelabel#1#2{\@bsphack
  \protected@write\@auxout{}%
         {\string\newlabel{#1.nme}{{#2}{#2}}}%
  \@esphack}

\makeatother

\begin{document}
\title{Generalized twisted quantum doubles of a finite group and rational orbifolds}

\author{Geoffrey Mason}
\address{Department of Mathematics, UC Santa Cruz, CA 95064}
\email{gem@ucsc.edu}
\author{Siu-Hung Ng}
\begin{quote}
Dedicated to Bob Griess on the occasion of his $71^{st}$ birthday.
\end{quote}
\thanks{The first author thanks the Simon Foundation, grant $\#427007$, for its support.}
\thanks{The second author was partially support by the NSF grant DMS 1501179.}
\address{Department of Mathematics, Louisiana State University, Baton Rouge, LA 70803}
 \email{rng@math.lsu.edu}
 \maketitle

 \begin{abstract} In previous work the authors introduced a new class of modular quasi-Hopf algebra
 $D^{\omega}(G, A)$, associated to a finite group $G$,  a central subgroup $A$ and a $3$-cocycle
 $\omega{\in}Z^3(G, \CC^{\times})$.\
 In the present paper we propose a description of the class of orbifold models of
 rational vertex operator algebras whose module category is tensor equivalent to $D^{\omega}(G, A)$-mod.\
 The paper includes background on quasi-Hopf algebras and a discussion of some relevant orbifolds.\\
 MSC Classification(2009){:}\ $16T99,18D99, 17B69$.
 \end{abstract}

\section{Introduction}
Since its introduction by Dijkgraaf, Pasquier and Roche \cite{DPR}, the \emph{twisted quantum double}
has been a source of inspiration in the related areas of \emph{quasiHopf algebras},
\emph{modular tensor categories} (MTC), and \emph{orbifold models} of holomorphic vertex operator algebras (VOA).
It follows from the work \cite{Mu1} of M\"{u}ger that the module category $D^{\omega}(G)$-mod of a twisted quantum double of a finite group
$G$ (notation and further details are provided below) is a MTC.\ Even before the idea
of a MTC existed, DPR had more-or-less conjectured \cite{DPR} that if $V$ is a holomorphic VOA admitting $G$ as a group of automorphisms then $V^G$-mod is equivalent to $D^{\omega}(G)$-mod for some $3$-cocycle
$\omega {\in} Z^3(G, \mathbb{C}^{\times})$.\ Kirillov considered the conjecture from the perspective of fusion categories \cite{Ki}. His work would imply the conjecture if one can show that all the $g$-twisted modules of $V$ form a $G$-graded fusion category. This condition amounts  to the rationality of $V^G$. With recent advances, this conjecture is now known for a large family
of groups $G$. The recent work \cite{MC} of Miyamoto and Carnahan proves that $V^G$ is  rational for any solvable group $G$. A complete solution to the conjecture seems to be within reach.

\medskip
On the other hand, much less known in the case of
\emph{rational orbifolds}, i.e., orbifolds $V^G$ where $V$ is a \emph{rational VOA}, but not necessarily holomorphic.\
Already in \cite{DPR}, the authors asked for a description of the $c{=}1$ ADE orbifolds
$V^G$ where $V{=}L(sl_2, 1)$ is the affine algebra of level $1$ associated to $sl_2$ (alternatively,
the rank $1$ lattice theory $V_{\sqrt{2}\mathbb{Z}}$ associated with the $A_1$ root lattice) and $G$ a finite subgroup
of $SO(3, \mathbb{R})$.\ Until recently there was no really satisfactory answer to this question.\
Fusion rules and $S$- and $T$-matrices for these theories have long been known (the icosahedral case proved to be particularly intractable) but a quasiHopf algebra replacing the
twisted double was missing.\ More generally, there does not seem to be even a conjectural description
of $V^G$-mod in the literature for any reasonably substantial class of rational orbifolds beyond the holomorphic case.

\medskip
In our recent work \cite{MN1} we introduced a \emph{generalization} of the twisted quantum double,
denoted by $D^{\omega}(G, A)$, which is a certain quasiHopf algebra quotient
of $D^{\omega}(G)$ obtained from a central subgroup $A$ of $G$.\ (The case $A{=}1$ reduces to $D^{\omega}(G)$.)\ We gave necessary and sufficient conditions that $D^{\omega}(G, A)$-mod
is a MTC.\ The purpose of the present paper is to present a conjectural description of those rational
orbifolds $V^G$ whose module category is equivalent to some $D^{\omega}(G, A)$-mod, and to discuss a few examples.\
These include the ADE examples mentioned above, thus providing an answer to the question raised by DPR.

\medskip
The paper is organized as follows. In Section 2 we present a general discussion of quasiHopf algebras,
with emphasis on the construction and properties of $D^{\omega}(G, A)$.\ In Section 3 we state our conjecture
relating $D^{\omega}(G, A)$ to certain rational orbifolds, and in Section 4 we consider some illustrative examples.

\medskip
We thank Terry Gannon for his interest and input, and for suggesting that this paper be written.

\section{QuasiHopf algebras and $D^{\omega}(G, A)$}
In this Section we provide some background, taken from \cite{MN1}, about the quasiHopf algebras $D^{\omega}(G, A)$, which we call \emph{generalized twisted quantum doubles}.

\medskip
We use the following notation for a finite group $G$.\
$\widehat{G}{=}\Hom(G, \mathbb{C}^{\times})$ is the group of \emph{characters} of $G$,
$x^g{=}g^{-1}xg\ (x, g{\in}G)$ is right conjugation in $G$, the \emph{centralizer} of $x$ in $G$ is
$C_G(x){:=}\{g{\in}G{\mid}x^g{=}x\}$, $Z(G){:=}\cap_x C_G(x)$ is the \emph{center} of $G$.

 \medskip
 We take our base field to be the complex numbers $\mathbb{C}$.\
 A \emph{quasiHopf algebra} is a tuple $(H, \Delta, \epsilon, \phi, \alpha, \beta, S)$, where
 $H$ is a unital algebra and $\Delta{:}H{\rightarrow} H\otimes H$ an algebra morphism that is
 \emph{quasicoassociative} in the sense that there is a map $\phi$
 (the \emph{Drinfeld associator}) making the following diagram commutative{:}
 \begin{eqnarray}\label{diag1}
 \xymatrix{
&H\otimes H \ar[ld]_{\Delta\otimes Id}\ar[rd]^{Id\otimes\Delta} \\
(H\otimes H)\otimes H\ar[rr]_{\phi} &&H\otimes(H\otimes H)
 }
\end{eqnarray}
 $S$ is the \emph{antipode}, $\epsilon$ the \emph{counit}, and $\alpha, \beta{\in}H$ are certain distinguished elements.\
 We generally suppress all reference to these elements of a quasiHopf algebra; this will not impair the reader's
 understanding of what  follows.\ One also requires $\phi$ to satisfy some \emph{cocycle} conditions in the form of certain diagrams involving fourfold tensor products of $H$ that are required to commute.\ Again we will generally suppress
 such details.\ A Hopf algebra is a quasiHopf algebra for which $\alpha{=}\beta{=}1$
 and $\phi{=}1\otimes 1\otimes 1$.\ For further background on quasiHopf algebras, see
 \cite{D}, \cite{K} and \cite{MN3}.

 \medskip
 One of the great virtues of a quasiHopf algebra $H$ is that the category $H$-mod of (finite-dimensional)
 $H$-modules is a finite tensor category, though not necessarily  a MTC.

 \medskip
 Fix a finite group $G$.\
The group algebra $\mathbb{C}G$ is a familiar example of a  Hopf algebra, the coproduct being defined by
$\Delta(g){=}g\otimes g\ (g{\in}G)$.\ Dualizing, the dual group algebra $\mathbb{C}^G$ has basis $e_g\ (g{\in}G)$
 \emph{dual} to the basis of group elements in $G$.\ It is a Hopf algebra with product and coproduct
 defined by
 \begin{eqnarray*}
&&e_ge_h{=}\delta_{g, h}e_g,\ \ \ \ \Delta{:}e_g \mapsto \sum_{ab{=}g} e_a \otimes e_b.
\end{eqnarray*}

 \medskip
 Now fix a normalized multiplicative $3$-cocycle $\omega{\in}Z^3(G, \mathbb{C}^{\times})$.\
A basic example of a quasiHopf algebra is the \emph{twisted dual group algebra} $\mathbb{C}^G_{\omega}$,
obtained from $\mathbb{C}^G$ just by replacing $1\otimes 1\otimes 1$
by a more interesting Drinfeld associator defined by multiplication by
\begin{eqnarray}\label{phidef}
\phi{:=}\sum_{a, b, c} \omega(a, b, c)^{-1}e_a\otimes e_b\otimes e_c.
\end{eqnarray}
Here, the cocycle conditions amount to the $3$-cocycle identity satisfied by $\omega$.

\medskip
The \emph{twisted quantum double}  $D^{\omega}(G){=}\mathbb{C}_{\omega}^G\otimes\mathbb{C}G$
occurs as the middle term of
a sequence of morphisms of quasiHopf algebras
 \begin{eqnarray*}
 \xymatrix{
\mathbb{C}_{\omega}^G\ar[rr]^i &&D^{\omega}(G)\ar[rr]^p&&\mathbb{C}G\\
 }
\end{eqnarray*}
where $i(e_g){=}e_g\otimes 1,\ p(e_g\otimes x){=}\delta_{g, 1}x$.\ $D^{\omega}(G)$ is itself a quasiHopf algebra
with the following product and coproduct{:}
\begin{eqnarray*}
&&(e_g\otimes x)\cdot (e_h\otimes y){=}\theta_g(x, y)\,\delta_{g^x, h}e_g\otimes xy\\
&&\Delta{:} e_g\otimes x \mapsto \sum_{ab=g} \gamma_x(a, b)\,e_a\otimes x \otimes e_b\otimes x.
\end{eqnarray*}
The scalars $\theta_g(x, y), \gamma_x(a, b)$ are determined by $\omega$ as follows{:}
\begin{eqnarray*}
&&\theta_g(x, y){:=} \frac{\omega(g, x, y)\omega(x, y, g^{xy})}{\omega(x, g^x, y)},\\
&&\gamma_x(a, b){:=} \frac{\omega(a, b, x)\omega(x, a^x, b^x)}{\omega(a, x , b^x)}.
\end{eqnarray*}

\medskip
We remark that the $2$-cochains defined by the $\theta$'s and $\gamma$'s have subtle properties which
govern much of the behaviour of $D^{\omega}(G)$ -- for example, the fact that $D^{\omega}(G)$
really is a quasiHopf algebra, which is not obvious.\ As another example, if we restrict $x, y$ to $C_G(g)$
then the 2-cochains $\theta_g, \gamma_g$ \emph{coincide} and become $2$-cocycles.\ We informally record this as
\begin{eqnarray}\label{theta=gamma}
\gamma_g{=}\theta_g{\in}Z^2(C_G(g), \mathbb{C}^{\times}),
\end{eqnarray}
where it is understood that the domains of $\gamma_g, \theta_g$ are here restricted to $C_G(g)$.

\medskip Since $D^{\omega}(G)$-mod is the Drinfeld center of $\Vec(G, \w)$, the fusion category of $G$-graded vector space with the associativity given by $\w$,
it follows from \cite{Mu1} that $D^{\omega}(G)$-mod is a MTC.

\medskip
We now fix another piece of data, namely a central subgroup $A{\subseteq}Z(G)$, and introduce
\begin{eqnarray*}
D^{\omega}(G, A){:=} \mathbb{C}_{\omega}^G\otimes\mathbb{C}(G/A).
\end{eqnarray*}
Notice that because the conjugation action of $A$ on $G$ is \emph{trivial} then the
\emph{identical} formulas used to define the operations in $D^{\omega}(G)$ still make sense
in $D^{\omega}(G, A)$.\ When $\w$ is \emph{compatible} with $A$, one can equip $D^{\omega}(G, A)$ with a product and coproduct similar to the case of $D^{\omega}(G)$ so that $D^{\omega}(G, A)$ is a quasiHopf algebra.

\medskip
Now it is natural to ask if, for suitable $\pi'$, there is a commuting  diagram of quasiHopf algebras and morphisms
 \begin{eqnarray}\label{pi'diag}
 \xymatrix{
&\mathbb{C}_{\omega}^G\ar[rr]^i\ar[d]_= &&D^{\omega}(G)\ar[rr]^p\ar[d]_{\pi'}&&\mathbb{C}G\ar[d]_{\pi}\\
&\mathbb{C}_{\omega}^G\ar[rr]^i &&D^{\omega}(G, A)\ar[rr]^{\bar{p}}&&\mathbb{C}(G/A)
 }
\end{eqnarray}
where $\pi(x){=}xA$ and $\bar{p}(e_g\otimes xA){=}\delta_{g, 1}xA$ for $x{\in}G$.
In case $\pi'$ exists, it will satisfy $\pi'(e_g\otimes x){=}\lambda e_g\otimes xA$ for a scalar $\lambda$ that
depends on $g$ and $x$.

\medskip
 Generally there will be no such $\pi'$, but we can give necessary and sufficient conditions for its existence.\
 To explain this we first consider the \emph{group-like elements} of $D^{\omega}(G)$.\ These are the (nonzero) elements
 $u{\in}D^{\omega}(G)$ such that $\Delta(u){=}u\otimes u$.\ As in the case of Hopf algebras,
 the group-like elements form a multiplicative subgroup $\Gamma^{\omega}(G){\subseteq} D^{\omega}(G)^{\times}$
 of the group of units.\ We are more interested in the \emph{central group-like elements},
 which is the subgroup $\Gamma_0^{\omega}(G){\subseteq} \Gamma^{\omega}(G)$
 consisting of elements that commute with all elements of $D^{\omega}(G)$.\ One can show (cf.\ \cite{MN3}) that there is a diagram of short exact sequences of groups (actually central extensions)
 \begin{eqnarray*}
 \xymatrix{
1\ar[r]& \widehat{G}\ar[rr]^i &&\Gamma^{\omega}(G)\ar[rr]^p&&B^{\omega}(G)\ar[r]&1\\
1\ar[r]& \widehat{G}\ar[rr]^i\ar[u]^{=} && \Gamma_0^{\omega}(G)\ar[u]\ar[rr]^p &&Z^{\omega}(G)\ar[r]\ar[u]&1
 }
\end{eqnarray*}
where
\begin{eqnarray*}
&&B^{\omega}(G){:=}\{g{\in}G{\mid}\gamma_g{\in}B^2(G, \mathbb{C}^{\times})\},\ \ \ Z^{\omega}(G){:=}B^{\omega}(G)\cap Z(G),
\end{eqnarray*}
and vertical arrows are containments.

\medskip
We emphasize that here, unlike the situation of (\ref{theta=gamma}), in order for
$g$ to lie in $B^{\omega}(G)$ it is necessary that the $2$-cochain $\gamma_g$ be a 2-coboundary on $G$
 rather than just the centralizer $C_G(g)$.\ On the other hand, if $g{\in}Z(G)$ then the context of (\ref{theta=gamma}) pertains,
so that $\gamma_g{=}\theta_g$ is always a $2$-cocycle on $G$, and the requirement to belong to
$Z^{\omega}(G)$ is that this 2-cocycle is in fact a 2-coboundary.

\medskip
It is shown in \cite{MN1} that the existence of the morphism $\pi'$ in (\ref{pi'diag}) is equivalent to
the existence of an enlarged diagram of central extensions
 \begin{eqnarray}\label{3ses}
 \xymatrix{
1\ar[r]& \widehat{G}\ar[rr]^i &&\Gamma^{\omega}(G)\ar[rr]^p&&B^{\omega}(G)\ar[r]&1\\
1\ar[r]& \widehat{G}\ar[rr]^i\ar[u]^{=} && \Gamma_0^{\omega}(G)\ar[u]\ar[rr]^p &&Z^{\omega}(G)\ar[r]\ar[u]&1\\
1\ar[r]&\widehat{G}\ar[rr]^i\ar[u]^{=}&& p^{-1}(A)\ar[u]\ar[rr]^p &&A\ar[r]\ar[u]&1\\
 }
\end{eqnarray}
where vertical maps are again containments and the \emph{lower ses splits}.\

\medskip
What is being asserted here is the following{:}\ the central subgroup $A{\subseteq}Z(G)$ is required to also lie in
$B^{\omega}(G)$, i.e., the $2$-cocycles $\theta_g$ for $g{\in}A$ are $2$-coboundaries on $G$.\ Moreover,
the ses obtained by pulling-back $A$ along $p$ must split.

\medskip
Once $\pi'$ is available, it follows that $D^{\omega}(G, A)$-mod is a subcategory of $D^{\omega}(G)$-mod.\
However, M\"{u}ger's theorem will generally not hold for $D^{\omega}(G, A)$, that is to say
$D^{\omega}(G, A)$-mod is generally \emph{not} a MTC.\ We will describe necessary and sufficient conditions,
established in \cite{MN1}, that make this so.

\medskip
First we say a bit more about the middle ses in (\ref{3ses}).\ Given $g{\in}Z^{\omega}(G)$ we have
$\gamma_g{=}\theta_g{\in}B^2(G, \mathbb{C}^{\times})$, so that there is $\tau_g{\in}C^1(G, \mathbb{C}^{\times})$
satisfying $\delta\tau_g{=}\theta_g$, i.e.,
\begin{eqnarray*}
\tau_g(x)\tau_g(y){=}\theta_g(x, y)\tau_g(xy)\ \ \ \ (x, y{\in}G).
\end{eqnarray*}
Because $\theta_g$ is a 2-coboundary, the twisted group algebra $\mathbb{C}^{\theta_g}G$ that it defines
is isomorphic to the group algebra $\mathbb{C}G$, and $\tau_g$ defines a choice of isomorphism
\begin{eqnarray*}
 \mathbb{C}G \stackrel{\cong}{\longrightarrow}\mathbb{C}^{\theta_g}G,\ \ \ x\mapsto \tau_g(x)x.
\end{eqnarray*}
There is no canonical choice of $\tau_g$, but any two of them differ by a character $\chi{\in}\widehat{G}$.\
Indeed, we have
\begin{eqnarray*}
\Gamma_0^{\omega}(G){=}
\left\{\sum_{x\in G} \tau_g(x)\chi(x)e_x\otimes g{\mid}\chi\in\widehat{G}, g\in Z^{\omega}(G)\right\}.
\end{eqnarray*}

\medskip
A 2-cocycle $\beta{\in}Z^2(Z^{\omega}(G), \widehat{G})$ that defines
the central extension that is the middle ses in (\ref{3ses}) is given by the formula
\begin{eqnarray*}
\beta(g, h)(k){=}\frac{\tau_g(k)\tau_h(k)}{\tau_{gh}(k)}\theta_k(g, h)\ \ \ \ \ (g, h{\in}Z^{\omega}(G), k{\in}G).
\end{eqnarray*}

\medskip
Because we are assuming that the lower ses in (\ref{3ses}) splits, the restriction of
$\beta$ to $A$ is a 2-coboundary on $A$.\ That is, there is a 1-cochain $\nu{\in}C^1(A, \widehat{G})$
such that
\begin{eqnarray*}
\beta(a, b){=}\frac{\nu(a)\nu(b)}{\nu(ab)}\ \ \ \ (a, b{\in}A).
\end{eqnarray*}

\medskip
Now we can show \cite{MN1} that the formula
\begin{eqnarray}\label{bichar}
(a{\mid}b)_{\nu}{:=}\frac{\tau_a(b)\tau_b(a)}{\nu(a)(b)\nu(b)(a)}
\end{eqnarray}
defines a \emph{symmetric bicharacter} $(\ {\mid}\ )_{\nu}{:}A{\times}A{\rightarrow} \mathbb{C}^{\times}$.\
Using results of M\"{u}ger \cite{Mu1}, \cite{Mu2}, it follows that
$D^{\omega}(G, A)$-mod is a MTC if, and only if, $(\ {\mid}\ )_{\nu}$ is \emph{nondegenerate}.\
Actually, $(\ {\mid}\ )_{\nu}$ is the restriction of a natural bicharacter defined on
$\Gamma_0^{\omega}(G)$, however we will not go into this here.


\section{Simple current orbifolds}
In the spirit of the proposal in  \cite{DPR} that the module categories of holomorphic orbifolds
coincide with the module categories $D^{\omega}(G)$-mod, in this Section we
describe those rational orbifolds $V^G$ expected to have the property that
$V^G$-mod is equivalent to some $D^{\omega}(G, A)$.\ In this setting, the case $A{=}1$
will reduce to the holomorphic orbifold case of DPR.

\medskip
We will generally be lapse about the detailed properties of the VOAs which we consider,
but we will be concerned with \emph{rational} VOAs $V$, which at the very least means that
$V$ is a \emph{simple} VOA and $V$-mod is a fusion category, i.e. a semisimple rigid tensor category with only finitely many simple objects.\
(See \cite{DLM1} for further background.)\ In fact we will only need a small subset of such theories, and to explain which ones these are we will review the theory of \emph{simple currents}.\

\medskip
A simple current is a simple (or irreducible) $V$-module, call it $M$, with the property that for all simple
$V$-modules $N$, the tensor product $M\boxtimes N$ is again a simple $V$-module.\ Otherwise stated,
$M$ represents an element in the Grothendieck group of $V$-mod, and an object of Frobenius-Perron dimension one.\
Thanks to the associativity of
$\boxtimes$, the distinct (isomorphism classes of) simple currents form a group with respect to tensor product
of modules, called the \emph{group of simple currents}.\ It is an abelian group because $V$-mod is also braided.\ The identity element is, of course, the vacuum space $V$.\ We are concerned here
with rational VOAs $V$ with the property that \emph{every simple $V$-module is a simple current} i.e., that $V$-mod is pointed.\
We call such a $V$ a \emph{simple current VOA}.\ There are many examples of such theories.\
In addition to holomorphic VOAs, where $V$ is the \emph{only} simple module, it is well-known
that \emph{lattice theories $V_L$} associated to a positive-definite even lattice $L$ are also simple current
VOAs.\ In this case, the group of simple currents is isomorphic to
$L^*/L$ where $L^*$ is the \emph{dual lattice} of $L$ \cite{D}.\

\medskip
Given a simple current VOA $V$, we may form the sum of \emph{all} simple $V$-modules $M$ to obtain a
larger space
\begin{eqnarray*}
\widetilde{V}{:=} \oplus_M M
\end{eqnarray*}
$\widetilde{V}$ can be equipped
with the structure of an \emph{abelian intertwining algebra} in the sense of Dong-Leowsky \cite{DL}.\
See \cite{DLM2} for further details.

\medskip\noindent
\underline{Conjecture}.\ Suppose that $V$ is a simple current VOA with group of simple currents $A$.\ Let $F{\subseteq}\Aut(V)$
be a finite group of automorphisms of $V$ (the \emph{orbifold group}) such that $F$ leaves invariant every simple $V$-module.\ Then
there is \emph{central extension}
 \begin{eqnarray*}
 \xymatrix{
1\ar[r]& A\ar[r] &G\ar[r]&F\ar[r]&1\\
 }
\end{eqnarray*}
and a $3$-cocycle $\omega{\in}Z^3(G, \mathbb{C}^{\times})$ such that $V^F$-mod is equivalent
to $D^{\omega}(G, A)$-mod as modular tensor categories.\ Conversely, if $D^{\omega}(G, A)$ exists and $D^{\omega}(G, A)$-mod
is a MTC, then there is a simple current VOA $V$ with group of simple currents $A$ and a group of automorphisms
$F{=}G/A$ such that $V^F$-mod is tensor equivalent to $D^{\omega}(G, A)$-mod.

\section{Examples}
We illustrate the Conjecture by discussing some examples in greater detail.\
Let the notation be as before.

\medskip\noindent
\underline{Example 1}.\ The holomorphic case.\\
Here, $A{=}1$ means that $V$ is a holomorphic VOA and $D^{\omega}(G, A){=}D^{\omega}(G)$.\ The Conjecture thus reduces that of DPR.

\medskip\noindent
\underline{Example 2}.\ The  case $F{=}1$, i.e., $G {=}A$ and $D^{\omega}(G, A){=}\mathbb{C}_{\omega}^G$ as quasiHopf algebras.\

\medskip
We saw before that in order for $\pi'$ to exist (i.e., $D^{\omega}(G, A)$ is a quasiHopf algebra quotient
of $D^{\omega}(G)$)  it is necessary that $A{\subseteq}Z^{\omega}(G)$,
meaning that each $\theta_g{\in}B^2(G, \mathbb{C}^{\times})$ is a $2$-coboundary for all $g {\in} A$.\ That is, $\omega{\in}Z^3(A, \mathbb{C}^{\times})_{ab}$ is an
\emph{abelian} $3$-cocycle on $A$ in the sense of \cite{MN3}, where such cocycles are studied extensively.\ They are closely related to
the abelian cohomology groups introduced by Eilenberg and Maclane \cite{EM}.\

\medskip
Thus in the case at hand,
the Conjecture asserts that $V$-mod is tensor equivalent to $\mathbb{C}^A_{\omega}$
for an abelian 3-cocycle $\omega$ that is nondegenerate in a suitable sense.\ Rather than explain what degeneracy means
here, we consider a special case in more detail.

\medskip
Let $L$ be an even lattice with bilinear form $\langle\ ,\  \rangle{:}L{\times}L{\rightarrow}\mathbb{Z}$, and let
$V_L$ be the corresponding lattice VOA.\ Then $A{=}L^*/L$.\ The conjecture says that
$V_L$-mod is tensor equivalent to the dual group algebra of $L^*/L$ twisted by an abelian 3-cocycle $\omega{\in}
Z^3(L^*/L, \mathbb{C}^{\times})$.\ The origin of $\omega$ is well-known in this case
(cf.\ \cite{DL}, Chapter 12 and \cite{MN3}, Section 11).\ Let $s{:}A{\rightarrow}L^*$ be a normalized section
of the canonical ses
 \begin{eqnarray*}
 \xymatrix{
1\ar[r]& L\ar[r] &L^*\ar[r]&A\ar[r]&1\\
 }
\end{eqnarray*}
Let $c_0{:}L^*{\times}L^*{\rightarrow}\mathbb{C}^{\times}$ be an \emph{alternating bicharacter}
with $c_0(\alpha, \beta){=}(-1)^{\langle \alpha, \beta\rangle}\ (\alpha, \beta{\in}L)$.\ Set
\begin{eqnarray*}
\omega(a, b, c){=}(-1)^{\langle s(c), s(a)+s(b)-s(a+b)\rangle}c_0(s(c), s(a){+}s(b){-}s(a+b)).
\end{eqnarray*}
$\omega$ is indeed an abelian 3-cocycle (cf.\ \cite{MN3}, Proposition 11.1) whose cohomology class is \emph{independent
of the choice of bicharacter and section}.

\medskip\noindent
\underline{Example 3}.\ The case $|A|{=}2$.\\
 Here, we are discussing simple VOAs $V$ with just two simple modules.\
One of them is the adjoint module $V$, the other we denote by $M$.\ Roughly speaking, we may think of
$\widetilde{V}{=}V{\oplus}M$ as a holomorphic \emph{super} VOA, though this may not  conform to some definitions.\ (This will not matter in the following discussion.)

\medskip
Next we discuss results of \cite{MN4} concerning the existence of $D^{\omega}(G, A)$ in the special case that $G$ is a finite group with
\emph{exactly one} subgroup of order\ $2$.\ We take $A$ to be the unique
subgroup of order $2$.\ Let $T$ be a $2$-Sylow subgroup of $G$.\ It is well-known
that $T$ is either cyclic or generalized quaternion.\ Groups with a unique involution have $2$-periodic cohomology by the Artin-Tate theory \cite{CE}, and for such groups the $2$-Sylow subgroup of $H^3(G, \mathbb{C}^{\times})$ is \emph{cyclic}
of order $|T|$.\ We call a $3$-cocycle $\omega{\in}Z^3(G, \mathbb{C}^{\times})$
a \emph{$2$-generator} if the corresponding cohomology class $[\omega]{\in}H^3(G, \mathbb{C}^{\times})$
has order \emph{divisible by $|T|$}.\

\medskip
In the setting of the previous paragraph, we can prove \cite{MN4} that $\pi'$ always exists, so that
$D^{\omega}(G, A)$ is a quasiHopf algebra quotient of $D^{\omega}(G)$.\ Moreover, $D^{\omega}(G, A)$-mod is a MTC if, and only if, $\omega$ is a $2$-generator.

\medskip
There are a number of interesting classes of finite groups with a unique subgroup $A$ of order $2$.\ The following lists some of them.
\begin{eqnarray}\label{gplist}
&&SL_2(q)\ \ (q\geq 3\ \mbox{an odd prime power}) \notag\\
&& 2.A_6, 2.A_7, 6.A_6, 6.A_7\\
&&\mbox{binary polyhedral groups}=\mbox{finite subgroups of}\ SU_2(\mathbb{C}) \notag
\end{eqnarray}
Here, $SL_2(q)$ is the group of $2{\times}2$ matrices of determinant 1 over the finite field of cardinality $q$, and
$2.A_n$ is the 2-fold perfect central extension of the alternating group
$A_n$.\ For $n{=}6, 7$ Schur discovered that there are exceptional 6-fold perfect central extensions
of $A_n$.\ There are a few overlaps among these groups{:}\ $SL_2(3)$ and $SL_2(5)$ are the binary
octahedral and icosahedral groups respectively, and $2.A_6\cong SL_2(9)$.\

\medskip
Our Conjecture says that for each pair $(G, \omega)$ such that $G$ has a unique subgroup of order $2$
and $[\omega]$ is a 2-generator of
$H^3(G, \mathbb{C}^{\times})$, there is a holomorphic super VOA $\widetilde{V}{=}V\oplus M$ (in the sense described before) such that $G{\subseteq}\Aut(\widetilde{V}), G/A{\subseteq}\Aut(V)$ and $V^{G/A}$-mod is tensor equivalent
to $D^{\omega}(G, A)$-mod.

\medskip
As far as we know, the existence of a $\widetilde{V}$ for most groups $G$ on the list (\ref{gplist}) is unknown.\
However, it is well-known that the binary polyhedral groups, and indeed $SU_2(\mathbb{C})$ itself, act on
the $c{=}1$ holomorphic super VOA defined by the $A_1$ root lattice theory $V_{\sqrt{2}\mathbb{Z}}$,
also known as the affine algebra $L_{sl_2, 1}$ of level 1 associated to the Lie algebra $sl_2(\mathbb{C})$.\
We review some of the details.

\medskip
Adopting the lattice perspective, the irreducible modules for $V_{\sqrt{2}\mathbb{Z}}$
consist of the adjoint module $V_{\sqrt{2}\mathbb{Z}}$ and the module $M{=}V_{1/\sqrt{2}+\sqrt{2}\mathbb{Z}}$
corresponding to the two cosets of $\sqrt{2}\mathbb{Z}$ in its dual lattice  $\frac{1}{\sqrt{2}} \mathbb{Z}$.\
Thus $\widetilde{V}{=}V_{\sqrt{2}\mathbb{Z}}{\oplus}V_{1/\sqrt{2}+\sqrt{2}\mathbb{Z}}{=}V_{\frac{1}{\sqrt{2}} \mathbb{Z}}$.\
We have $\Aut(V){=}SO_3(\mathbb{R})$, obtained by exponentiating the weight $1$ states of $V$ (which
form the Lie algebra $sl_2$).\ The projective action of this group on $M$ lifts to a linear action of its
universal cover $SU_2(\mathbb{C})$, which is the full automorphism group of $\widetilde{V}$.

\medskip
The binary polyhedral groups $G$ of even order are the finite subgroups of $SU_2(\mathbb{C})$ that contain
the center $\{\pm 1\}$.\ Setting $A{=}\{\pm 1\}$, the preceding discussion shows that
$G{\subseteq}\Aut(\widetilde{V})$ and $G/A{\subseteq}\Aut(V)$, so we have rational orbifolds
$V^{G/A}$.\ These theories have been studied extensively in both the physical and mathematical literature,
and the $S$- and $T$-matrices are known.\ See, for example, \cite{CG} and \cite{DV} for type $D$.\
Dong and Nagatomo \cite{DN} computed the fusion rules  for $V^{G/A}$ on the basis of the VOA axioms in all cases except the icosahedral example.

\medskip
Now we can compute the  $S$- and $T$-matrices for all
$D^{\omega}(G, A)$ whenever $G$ is a binary polyhedral group and
$\omega$ is a $2$-generator.\ Note that $H^3(G, \mathbb{C}^{\times})$ is cyclic
of order $|G|$ in this case, so that the number of $2$-generators is exactly $|G|/2$.\
Consider, for example the icosahedral group $G{=}SL_2(5)$.\ This is
a group of order 120, so there are $60$
$2$-generators $[\omega]$ of order $8, 24, 40$ or $120$.\ In the case of $V^{G/A}$-mod discussed
above, the $S$-matrix appears to correspond to an $[\omega]$ of maximal order $120$.

\medskip
 In this way we get a large number of
MTCs and accompanying modular data, and the Conjecture says that among them we should find
 the categories $V^{G/A}$-mod.\ This appears to indeed be the case, since in particular we can find matching $S$-matrices in all cases.\

\end{document}